\documentclass[12pt,fleqn]{article}

\usepackage[english]{babel}
\usepackage{makeidx}
\usepackage{latexsym,amsfonts,amssymb,amsmath,longtable,amsthm}
\input amssym.def
\usepackage{latexsym}
\usepackage[all]{xy}
\usepackage{amsfonts}
\usepackage{amsmath}
\usepackage{mathrsfs}
\usepackage{color}
\usepackage{ulem}
\usepackage[table]{xcolor}
\usepackage[unicode,breaklinks=true,colorlinks=true]{hyperref}

\catcode`@=11
\newbox\tr@tto
\setbox\tr@tto=\hbox{{\count0=0\dimen0=-
,9pt\dimen1=1,1pt\loop\ifnum\count0<11 \advance \count0 by1 \vrule
width.51pt height\dimen1
     depth\dimen0\kern-0.17pt\advance\dimen0 by-
0.05pt\advance\dimen1
     by0.1pt\repeat \loop\ifnum\count0<21\advance \count0 by1
\vrule
     width.6pt height\dimen1 depth\dimen0\kern-0.2pt
\advance\dimen0
     by-0.1pt\advance\dimen1 by 0.05pt\repeat}}
\def\medint{\displaystyle\copy\tr@tto\kern-10.4pt\int}
\catcode`@=12

\def\Xint#1{\mathchoice
   {\XXint\displaystyle\textstyle{#1}}%
   {\XXint\textstyle\scriptstyle{#1}}%
   {\XXint\scriptstyle\scriptscriptstyle{#1}}%
   {\XXint\scriptscriptstyle\scriptscriptstyle{#1}}%
   \!\int}
\def\XXint#1#2#3{{\setbox0=\hbox{$#1{#2#3}{\int}$}
     \vcenter{\hbox{$#2#3$}}\kern-.5\wd0}}

\def\dashint{\Xint-}

\newcommand{\R}{{\mathbb R}}

\newcommand{\N}{{\mathbb N}}
\newcommand{\Z}{{\mathbb Z}}

\newcommand{\al}{{\alpha}}
\newcommand{\Le}{{\mathscr L}}

\newcommand{\M}{{\mathcal  M}}

\renewcommand{\H}{{\mathcal H}}
\newcommand{\rank}{{\mathrm{rank\,}}}

\newcommand{\p}{{p_\circ}}

\renewcommand{\t}{\tau}
\renewcommand{\tt}{\tau_*}

\newcommand{\loc}{{\rm loc}}

\newcommand{\meas}{\mathop{\rm meas}}
\newcommand{\diam}{\mathop{\rm diam}}

\newcommand{\LL}{\mathrm{L}}
\newcommand{\WW}{\mathrm{W}}

\newcommand{\CC}{\mathrm{C}}
\newcommand{\cc}{C^{k,\alpha+}}
\newcommand{\Cc}{C^{k,\alpha}}
\newcommand{\ZZ}{\widetilde Z}

\newcommand{\Lll}{\Le^{k+\al}_{p,1}}
\newcommand{\LK}{\Le^{\al}_{p}}

\newtheorem{ttt}{\bf Theorem}[section]

\newtheorem{cor}{\bf Corollary}[section]

\theoremstyle{remark}
\newtheorem{rem}{\bf Remark}[section]

\numberwithin{equation}{section}

\newcommand{\dd}{{\rm d}}

\title{Morse--Sard theorem and Luzin $N$-property: a~new synthesis result for Sobolev spaces}
\author{Adele Ferone, \, Mikhail V.~Korobkov,  \, and \, Alba Roviello}

\begin{document}

\maketitle

\begin{abstract}
For a regular (in a sense) mapping $v:\R^n\to\R^d$ we study the
following problem:\  {\sl let $S$ be a subset of $m$-critical a
set~$\ZZ_{v,m}=\{\rank \nabla v\le m\}$ and the equality
$\H^\tau(S)=0$
 (or the~ inequality~ $\H^\tau(S)<\infty$\,) holds for some
 $\t>0$. Does it imply that $\H^{\sigma}(v(S))=0$ for some $\sigma=\sigma(\tau,m)$?}
(Here $\H^\t$ means the $\t$-dimensional Hausdorff measure.)

For the classical classes $C^k$-smooth and $C^{k+\al}$-Holder
mappings this problem was solved in the papers by Bates and Moreira. We solve
the problem for Sobolev $W^k_p$ and fractional Sobolev
$W^{k+\al}_p$ classes as well. Note that we study the Sobolev case
under minimal integrability assumptions $p=\max(1,n/k)$, i.e., it
guarantees in general only {\it the continuity} (not everywhere
differentiability) of a~mapping.

In particular, there is an interesting and unexpected analytical
phenomena here: if $\t=n$ (i.e., in the case of Morse--Sard
theorem), then the value $\sigma(\tau)$ is the same for
the~Sobolev $W^k_p$ and for the classical $C^k$-smooth case. But
if $\t<n$, then the value $\sigma$ depends on~$p$ also; the
value~$\sigma$ for $C^k$ case could be obtained as the limit
when~$p\to\infty$. The similar phenomena holds for Holder
continuous $C^{k+\al}$ and for the fractional Sobolev
$W^{k+\al}_p$ classes.

The proofs of the most results are based on our previous joint
papers with J.~Bourgain and J.~Kristensen  (2013, 2015). We also
crucially use very deep Y.~Yomdin's entropy estimates of near
critical values for polynomials (based on algebraic geometry
tools).
\medskip



\end{abstract}


\section{Introduction}\label{Introd}

The Morse--Sard theorem in its classical form states that the
image of the set of critical points of a $\CC^{n-d+1}$ smooth
mapping $v\colon \R^n \to \R^d$ has zero Lebesgue measure in
$\R^d$. More precisely, assuming that $n \geq d$, the set of
critical points for $v$ is $Z_v=\{x\in\R^n\,:\,\rank \nabla
v(x)<d\}$ and the conclusion is that
\begin{equation}\label{classical}
\Le^d(v(Z_v))=0
\end{equation}
whenever $v\in C^k$ with $k\ge\max(1,d-m+1)$. The theorem was proved by Morse~\cite{Mo} in 1939 for the case
$d=1$ and subsequently by Sard~\cite{S} in 1942 for the general
vector--valued case. The celebrated results of Whitney~\cite{Wh}
show that the $\CC^{n-d+1}$ smoothness assumption on the mapping
$v$ is sharp.

Another important item of the real analysis, $N$-property, means
that the image $v(E)$ has zero measure whenever $E$ has zero
measure (see the recent paper~\cite{FKR-n}, where we discuss the
history of the~topic).

We need some usual notation. Fix a pair of positive parameters
$\tau$ and $\sigma$.   A continuous mapping  $v:\R^n\to\R^d$   \
is said to satisfy \,$(\tau,\sigma)$-$N$-{\it property,} if
$$\H^\sigma(v(E))=0\mbox{ \ whenever \ $\H^\tau(E)=0$,}$$ where
$\H^\tau$ means the Hausdorff measure.

For  a $C^1$-smooth mapping $v:\R^n\to\R^d$  and for an integer
number $m\in\Z_+$ denote
$$\ZZ_{v,m}=\{x\in\R^n\,:\,\rank \nabla v(x)\le m\}.\footnote{We
use the symbol $\ZZ$, since in our previous papers we denoted
$Z_{v,m}=\{x\in\R^n\,:\,\rank \nabla v(x)<m\}$. So in the present
notation \,$\ZZ_{v,m}=Z_{v,m+1}$.}$$ Then for parameters
$\tau,\sigma>0$ we say that that a~mapping $v:\R^n\to\R^d$ \
satisfies $(\tau,\sigma,m)$-$N$-{\it property,} if
$$\H^\sigma(v(E))=0\mbox{ \ whenever \ $E\subset \ZZ_{v,m}$ \ with
\ $\H^\tau(E)=0$.}$$

Further, we say that that a~mapping $v:\R^n\to\R^d$ \ satisfies
{\it strict $(\tau,\sigma,m)$-$N$-property,} if
$$\H^\sigma(v(E))=0\mbox{ \ whenever \ $E\subset \ZZ_{v,m}$ \ with
\ $\H^\tau(E)<\infty$.}$$

Using this notation, the above classical Morse--Sard theorem
means, that every $C^k$-mapping $v:\R^n\to\R^d$ has strict
$(n,d,d-1)$-$N$-property if $k\ge n-d+1$.

The starting point for our research is the following recent result for classically smooth case.

\begin{ttt}[Bates S.M. and Moreira C., 2002  \cite{BM,CGTAM}]\label{MSN}{\sl Let
$m\in\{0,\dots,n-1\}$, \,$k\ge1$, \,$d\ge m$, \,$0\le\al\le1$,
\,and \,$v\in \Cc(\R^n,\R^d)$. Then for any $\tau\in [m,n]$ the
mapping $v$ has $(\tau,\sigma,m)$-$N$-property with
 \begin{equation}\label{htf1}
\sigma=m+\frac{\tau-m}{k+\alpha}.
\end{equation}
Moreover, this $N$-property is strict  if at least one of the
following additional assumptions is fulfilled:

1)  $\tau=n$ (in particular, it includes  the case of the
classical Morse--Sard theorem);

2) $\t>m$ and $\alpha=0$ (that means $v\in C^k$);

3) $\t>m$ and $v\in C^{k,\alpha+}(\R^n,\R^d)$.}
\end{ttt}

Here we  say that a~mapping $v:\R^n\to\R^d$ belongs to the class
$\Cc$ for some positive integer $k$ and $0<\alpha\le1$ \  if \
$v\in C^k$ \,and  there exists a~constant $L\geq 0$ such that
\begin{equation*}\label{hic1}\mbox{$|\nabla^kv(x)-\nabla^kv(y)|\le
L\,|x-y|^\alpha$ \qquad for all $x,y\in\R^n$.}
\end{equation*}To simplify the notation, let us make the following agreement: for
$\alpha=0$ we identify~$C^{k,\alpha}$ with usual spaces of
$C^k$-smooth mappings.

Analogously, we say that
 a~mapping $v:\R^n\to\R^d$ belongs to the class $\cc$
for some positive integer $k$ and $0<\alpha\le1$, if there exists
a function~$\omega:\R_+\to\R_+$ such that $\mbox{$\omega(r)\to0$ \
as \,$r\to0$}$ \ and \begin{equation}\label{uuub}
\mbox{$|\nabla^kv(x)-\nabla^kv(y)|\le\omega(r)\cdot |x-y|^\alpha$
\ \ \quad whenever \ $|x-y|<r$}.
\end{equation}

Note that the assertion of Theorem~\ref{MSN} is rather sharp: for example, if its conditions 1)--3) are not satisfied, than 
the corresponding~$(\tau,\sigma,m)$-$N$-property is not strict in general, it follows from  Whitney's counterexamples~\cite{Wh}, see also~\cite{Nor} for commentaries.

Of course, the assertion of Theorem~\ref{MSN} includes Morse--Sard theorem and many
other results on this topic as partial cases; for convenience, we
made some historical references below in Subsection~\ref{Hist}.
The purpose of our paper is to extend this result to the mappings
of Sobolev spaces.

\subsection{Morse--Sard--Luzin type theorem for the case of Sobolev spaces}\label{bst0}

In this subsection  $\WW^k_{p}(\R^n,\R^d)$ means the space of
Sobolev mappings with all derivatives of order $j\le k$ belonging
to the Lebesgue space~$\LL_{p}$.

Let $k\in\N$, $1<p<\infty$  and $0\le\alpha<1$. One of the most
natural type of
 fractional Sobolev spaces  is
{\it (Bessel) potential spaces} $\Le^{k+\alpha}_p$.

Recall, that a~function $v:\R^n\to\R^d$ belongs to the space
$\Le^{k+\alpha}_p$, if it is a convolution of a function~$g\in
L_p(\R^n)$ with the Bessel kernel~$G_{k+\alpha}$, where
$\widehat{G_{k+\alpha}}(\xi)=(1+4\pi^2\xi^2)^{-(k+\alpha)/2}$. It
is well known that for the integer exponents (i.e., when $\al=0$)
one has the identity
\begin{equation*}\label{fN1} \Le^{k}_p(\R^n)=W^k_p(\R^n)\qquad\mbox{ if }\ \quad 1<p<\infty.
\end{equation*}
As well-known, if   $(k+\al)p>n$, then functions from the
potential space $\Le^{k+\alpha}_p(\R^n)$ are continuous by Sobolev
Imbedding theorem, but in general the gradient~$\nabla v$ is not
well-defined everywhere. Thus now  for the Sobolev case the
$m$-critical set is defined as
$$\ZZ_{v,m}=\{x\in\R^n:x\in A_v\mbox{\ \ or \ }x\in\R^n\setminus A_v\  \ \mbox{with}\,\ \rank\nabla v(x)\le m\}.$$
Here $A_v$ means  the set of `bad' points at which either the
function~$v$ is not differentiable or which are not the  Lebesgue
points for~$\nabla v$. So in the paper\footnote{In our previous
papers we consider the $m$-critical points and `bad'
points $A_v$ separately.}
 we consider these `bad' nonregular points automatically as $m$-critical for any~$m$ (such assumption, of course, makes the corresponding $(\tau,\sigma,m)$-$N$-properties more stronger).

\begin{ttt}\label{MSN-S}{\sl Let
$m\in\{0,\dots,n-1\}$, \,$k\ge1$, \,$d\ge m$, \,$0\le\al<1$,
\,$p>1$, \,$(k+\al)p>n$, \,and let \,$v\in
\Le^{k+\al}_{p}(\R^n,\R^d)$. Denote~$\tt=n-(k+\al-1)p$. Suppose in
addition that
\begin{equation*}\label{st1} \t>m\qquad\mbox{ and }\qquad\t>\tt,
\end{equation*}
then the mapping~$v$  has {\bf strict}
$(\tau,\sigma,m)$-$N$-property with
\begin{equation}\label{stf1} \sigma=m+\frac{p(\t-m)}{\t+(k+\al)p-n}.
\end{equation}
Further, if $\t=m>\tt$, then $v$ has nonstrict
$(\tau,m,m)$-$N$-property. }
\end{ttt}

We need to make several remarks here.
\begin{itemize}
\item  First of all, let us note, that the value~$\sigma$ in
Theorems~\ref{MSN}--\ref{MSN-S} coincide for the boundary cases
$\t=m$ or $\t=n$, but they are different for $m<\t<n$ (of course,
then $\sigma$ for Sobolev case is larger). Nevertheless, $\sigma$
in Theorem~\ref{MSN} could be obtained by taking a limit
in~(\ref{stf1}) as
 $p\to\infty$;

 \item
Recall, that by approximation results (see, e.g., \cite{Sw} \,and
\,\cite{KK15}\,) the set of `bad' points~$A_v$ is rather small,
i.e., it has the Hausdorff dimension~$\tt$:
\begin{equation}\label{small1}\H^\tau(A_v)=0\qquad\forall\tau>\tt:=n-(k+\al-1)p\qquad\mbox{ if \ $v\in\Le^{k+\al}_p(\R^n)$}.
\end{equation}
In particular, $A_v=\emptyset$ if $(k+\al-1)p>n$.

\item The condition~$\t>\tt$ in Theorem~\ref{MSN-S} is essential
and sharp: namely,  in the paper~\cite{FKR-n} we constructed a
counterexample of a mapping from $\Le^{k+\al}_p(\R^n)$ not
satisfying the $(\tau,\sigma,m)$-$N$-property with
$\t=\tt=m=\sigma=1$.

\item The usual $(\tau,\sigma)$-$N$-properties (without
constraints on the gradient, i.e., when $m=n$) were studied in our
previous paper~\cite{FKR-n}, see also subsection~\ref{Hist}, \,Theorems~\ref{LPT1}--\ref{LPT2}. (One
has to use these usual~$N$-properties also if the
assumptions~$\t>m$ and $\t>\tt$ of Theorem~\ref{MSN-S} are not
satisfied.)
\end{itemize}

Thus above Theorem~\ref{MSN-S} omits the limiting cases
$(k+\al)p=n$ and $\tau=\tt$. However, it is possible to cover
these cases as well using the Lorentz norms. Namely, denote by
$\Le^{k+\alpha}_{p,1}(\R^n,\R^d)$  the space of functions which
could be represented as a~convolution of the~Bessel
potential~$G_{k+\al}$ with a function~$g$ from the Lorentz
space~$L_{p,1}$ (see the definition of these spaces in the
section~\ref{prel}).

\begin{ttt}
\label{MSN-SL}{\sl  Let $m\in\{0,\dots,n-1\}$, \,$k\ge1$, \,$d\ge
m$, \,$0\le\al<1$, \,$p\ge1$ \,and \,let $v:\R^n\to \R^d$ be
a~mapping for which one of the following cases holds:
\begin{itemize}
\item[(i)] \,$\al=0$, \,$k\ge n$, \,and \,$v\in W^k_1(\R^n,\R^d)$;

\item[(ii)] \,$0\le\al<1$, \,$p>1$, \,$(k+\al)p\ge n$, \,and
\,$v\in \Le^{k+\al}_{p,1}(\R^n,\R^d)$.
\end{itemize}
Denote $\tt=n-(k+\al-1)p$. Suppose in addition that
\begin{equation*}\label{st1} \t>m\qquad\mbox{ and }\qquad\t\ge \tt,
\end{equation*}
then the mapping~$v$  has {\bf strict}
$(\tau,\sigma,m)$-$N$-property with the same~$\sigma$ defined
by~(\ref{stf1}). Further, if $\t=m\ge\tt$, then $v$ has the
corresponding  nonstrict $(\tau,m,m)$-$N$-property. }
\end{ttt}

So here the limiting case $\t=\tt$ is included. Some other
commentaries:
\begin{itemize}
\item Recall, that by approximation results (see, e.g., \cite{Sw}
\,and \,\cite{KK15}\,) the set of `bad' points $A_v$ for this
Sobolev--Lorentz case has the same Hausdorff
dimension~$\tt=n-(k+\al-1)p$, but it is smaller in a sense,
namely:
\begin{equation}\label{small2}\H^{\tt}(A_v)=0\qquad\mbox{ if \ $v$ is from Theorem~\ref{MSN-SL} }.
\end{equation}
(compare with~(\ref{small1})\,). In particular, $A_v=\emptyset$ if
$(k+\al-1)p\ge n$.

\item For  the integer exponents (i.e., when $\al=0$) the
Sobolev--Lorentz potential space has a~more simple and natural
description:
\begin{equation*}\label{fN1} \Le^{k}_{p,1}(\R^n)=W^k_{p,1}(\R^n)\qquad\mbox{ if }\ \quad 1<p<\infty,
\end{equation*} there by  $\WW^k_{p,1}$ we denote
the subspace of $W^k_{p}$ consisting of functions whose
derivatives of order~$k$ belongs to the Lorentz space $L_{p,1}$
(see, e.g.,~\cite{FKR-n}).
\end{itemize}

\subsection{Some historical remarks}
\label{Hist}

There are a lot of papers devoted to the Morse--Sard theorem, and
the above formulated results includes many previous theorems as
partial cases. For example, for smooth case if $\alpha=0$, $\t=n$, then we have
$$\sigma=m+\frac{n-m}{k},$$
and the assertion of Theorem~\ref{MSN} coincides with the
classical Federer--Dubovitski\u{\i} theorem, obtained
 almost simultaneously by
Dubovitski\u{\i}~\cite{Du2} in 1967 and Federer \cite[Theorem
3.4.3]{Fed} in 1969.  Of course, it includes the original
Morse--Sard theorem as partial case (when $k=n-m,\sigma=m+1$\,). 

Note also, that Theorem~\ref{MSN} was formulated as a Conjecture by A.Norton in~\cite[page 369]{Nor} and it  includes as partial cases some relative results of other mathematicians: 
Norton himself (who proved the assertion for the case~$\sigma=d$, $\t=(k+\al)(d-m)+m$\,), Y.Yomdin~\cite{Yom} (case~$\tau=n$, \,$v\in C^{k,\al+}$, see also~\cite{BHS}\,),
M.~Kucera~\cite{Kucera} (case $\t=n$ , $m=1$, i.e., when the gradient totally vanishes on the critical set), etc.

\

Concerning the Sobolev case, in the pioneering paper by De~Pascale~\cite{DeP} the assertion of
the initial Morse--Sard theorem~(\ref{classical}) \ (i.e., when
$k=n-d+1$, $m=d-1$, $\sigma=d$\,) was obtained for the Sobolev
classes $W^k_p(\R^n,\R^m)$ under additional assumption $p>n$ (in
this case the classical embedding $W^k_p(\R^n,\R^m)\hookrightarrow
C^{k-1}$ holds, so there are no problems with nondifferentiability
points).

 Some other Morse--Sard type theorems for Sobolev cases were
obtained in~\cite{BHS} and ~\cite{HZ}, these papers mainly concern
the Dubovitski\u{\i}--Fubini type properties for the Morse--Sard
theorem, which will be discussed  in the next subsection.

In addition to the above mentioned papers there is a growing
number of papers on the topic, including
\cite{AZ,AS,Barbet,Bates,Bu,Fig,PZ,Putten,Putten1}.

Finally, Theorems~\ref{MSN-S} and \ref{MSN-SL} for the most
important case $\tau=n$ were obtained in our previous
paper~\cite{FKR-MS} (see also our preceding
articles~\cite{BKK,BKK2,HKK,KK3,KK15} of the first author with J.Bourgain, J.Kristensen,  and P. Haj\l{}asz on this topic).

The usual $(\tau,\sigma)$-$N$-properties (without constraints on
the gradient, i.e., when $m=n$) were studied in our previous
paper~\cite{FKR-n}, where we proved the following two theorems:

\begin{ttt}[\cite{FKR-n}]
\label{LPT1}{\sl Let $\al>0$, \,$1<p<\infty$, \,$\alpha p>n$,  and
$v\in \LK(\R^n,\R^d)$. Suppose that $0<\tau\le n$. Then   the
following assertions hold:
\begin{itemize}
\item[($\circ$)] \,if $\tau\ne\tau_*=n-(\alpha-1)p$, then $v$ has
the~$(\tau,\sigma)$-$N$-property,  where the value
$\sigma=\sigma(\tau)$ is defined as
\begin{equation}\label{np1}
\sigma(\tau):=\left\{\begin{array}{lcr}\tau, \ \ &\mbox{\rm if \
}& \tau\ge\tau_*:=n-(\alpha-1) p;\\ [13pt] \frac{p\,\tau}{\alpha
p-n+\tau}, \ \ &\mbox{\rm if \ }& 0<\tau
<\tau_*.\end{array}\right.
\end{equation}

\item[($\circ$$\circ$)] \,if $\al>1$ \,and \,$\tau=\tau_*>0$ then
$\sigma(\tau)=\tt$ and  the mapping $v$ in general has no \
$(\tt,\tt)$-$N$-property, i.e., it could be $\H^{{\tt}}(v(E))>0$
for some $E\subset\R^n$ with $\H^{{\tt}}(E)=0$.
\end{itemize}
}
\end{ttt}

The similar results were announced in~\cite{Alb}, see \cite{FKR-n}
for our commentaries and other historical remarks on this
important case.

The above Theorem~\ref{LPT1} omits the limiting cases $\alpha p=n$
and $\tau=\tt$. As above, it is  possible to cover these cases as
well using the Lorentz norms.

\begin{ttt}[\cite{HKK,FKR-n}]
\label{LPT2}{\sl Let  $v:\R^n\to \R^d$ be
a~mapping for which one of the following cases holds:
\begin{itemize}
\item[(i)] \,$v\in W^k_1(\R^n,\R^d)$ for some $k\in\N$, $k\ge n$;

\item[(ii)] \, $v\in \Le^{\al}_{p,1}(\R^n,\R^d)$ \ for  some $\al>0$, \ $p\in (1,\infty)$ with $\al p\ge n$.
\end{itemize}
Suppose that $0<\tau\le n$. Then $v$
is a~continuous function satisfying
the~$(\tau,\sigma)$-$N$-property,  where again the value
$\sigma=\sigma(\tau)$ is defined in~(\ref{np1}) \,(with  $\al=k$ and $p=1$ for the~(i) case). }
\end{ttt}

So, in the last theorem the critical case $\tau=\tt$ is {\bf included}.

\subsection{The Dubovitski\u{\i}--Fubini type properties for the Morse--Sard
theorem}\label{s-Dub}

As it was mentioned by A.Norton~\cite[page 369]{Nor}, the absence of a Fubini theorem
for Hausdorff measure makes an obstacle for proofs  of some new
Morse--Sard type theorems.  Nevertheless, in 1957 Dubovitski\u{\i}
proved, that surprisingly some Fubini type properties always hold
for the Morse--Sard topic.

\vspace{2mm}

 \noindent {\bf Theorem A
(Dubovitski\u{\i}~1957~\cite{Du})}.\label{AA} {\sl Let $n,d,k \in
\N$, and let
 $v \colon \R^n\to\R^d$ be a
$\CC^k$--smooth mapping. Then
\begin{equation}\label{dub1}
\H^\mu (Z_{v}\cap v^{-1}(y))=0\qquad\mbox{ for $\Le^d$-a.a.
}y\in\R^d,
\end{equation}
where $\mu=n-d+1-k$ and $Z_v=\{x\in\R^n:\rank\nabla v(x)<d\}$.
}\vspace{2mm}

Here and in the following we interpret $\H^\beta$ as the counting
measure when $\beta\le0$. Thus for $k\ge n-d+1$ we have $\nu
\le0$, and $\H^\mu$ in (\ref{dub1}) becomes simply the counting
measure, so the Dubovitski\u{\i} theorem contains the Morse--Sard
theorem as particular case.

It turns out that the similar Fubini type extensions hold for the
Theorems~\ref{MSN}--\ref{MSN-SL} stated above.

\begin{rem}\label{Complex}
The following language below may seem too technical and cumbersome. So,
a disinterested reader can omit them; anyway  the main results of
the article are the above theorems~\ref{MSN-S}--\ref{MSN-SL}.
Nevertheless, authors consider the following theorems as important
strengthens of theorems~\ref{MSN}--\ref{MSN-SL}, as they allow to
realise the idea of Dubovitsky's approach  in general situation,
and include all the theorems given in this article as a
 particular case; moreover, they are new even for the classical smooth cases $C^k$ and $C^{k,\al}$. 
\end{rem}

We need some notation. For parameters $\mu\ge0$, $q\ge m$,
$\tau>0$ we say that that a~mapping $v:\R^n\to\R^d$ \ satisfies
$(\tau,\mu,q,m)$-$N$-{\it property,} if
\begin{equation}\label{def-dub}
\H^\mu(E\cap v^{-1}(y))=0 \mbox{ \ for \,$\H^q$-almost all $y\in
v(E)$ \ whenever \ $E\subset \ZZ_{v,m}$ \ with \
$\H^\tau(E)=0$.}\end{equation}
Recall, that here as above
$\ZZ_{v,m}=\{x\in\R^n\,:\,\rank \nabla v(x)\le m\}.$ Obviously,
\begin{equation}\label{dub3}
\mbox{if $\mu\le0$, \ then the $(\tau,\mu,q,m)$-$N$-property is
equivalent to the $(\tau,q,m)$-$N$-property.}
\end{equation}

Further, we say that that a~mapping $v:\R^n\to\R^d$ \ satisfies
{\it strict $(\tau,\mu,q,m)$-$N$-property,} if
$$\H^\mu(E\cap v^{-1}(y))=0 \mbox{ \ for \,$\H^q$-almost all $y\in
v(E)$ \ whenever \ $E\subset \ZZ_{v,m}$ \ with \
$\H^\tau(E)<\infty$.}$$

\begin{ttt}[Smooth case $v\in C^{k,\al}(\R^n,\R^d)$]\label{MSN-D}{\sl
\ \ Under assumptions of Theorem~\ref{MSN} one can replace the
assertion about $(\t,\sigma,m)$-$N$-properties
 by the more strong assertion about
$(\tau,q,\mu,m)$-$N$-property for any $\t\in[m,n]$ and
$q\in[m,\sigma]$ with
\begin{equation}\label{h-d-f} \mu=\t-m-(k+\al)(q-m).
\end{equation}
Further, if $q>m$ and at least one of the corresponding conditions
{\rm1)--3)} of Theorem~\ref{MSN} is fulfilled, then this
$(\tau,q,\mu,m)$-$N$-property is strict.}
\end{ttt}

The similar assertions hold for Sobolev and Sobolev--Lorentz cases
(we use the definition from subsection~\ref{bst0} for the
$m$-critical set $\ZZ_{v,m}$ of Sobolev functions).

\begin{ttt}[Sobolev case $v\in \Le^{k,\al}_p(\R^n,\R^d)$, $(k+\al)p>n$]\label{MSN-S-D}{\sl
\ Under assumptions of Theorem~\ref{MSN-S} one can replace the
assertion about strict $(\t,\sigma,m)$-$N$-properties
 by the more strong assertion about
strict $(\tau,q,\mu,m)$-$N$-property for any $\t>\max(\tt,m)$,
\,$q\in(m,\sigma]$  with
\begin{equation}\label{stf1-D} \mu=\t-m-(k+\al-\frac{n}p+\frac{\tau}p)(q-m).
\end{equation}
Further, if $q=m$, $\t>\tt$, and $\t\ge m$,   then $v$ has
nonstrict $(\tau,m,\mu,m)$-$N$-property with $\mu=\tau-m$. }
\end{ttt}

\begin{ttt}[Sobolev--Lorentz case $v\in \Le^{k,\al}_{p,1}(\R^n,\R^d)$, $kp\ge
n$]\label{MSN-SL-D}{\sl \
Under assumptions of Theorem~\ref{MSN-SL} one can replace the
assertion about strict $(\t,\sigma,m)$-$N$-properties
 by the more strong assertion about strict
$(\tau,q,\mu,m)$-$N$-property for any $\t\ge\tt$, $\t>m$,
\,$q\in(m,\sigma]$, and with the same~$\mu$ as in~(\ref{stf1-D}).
Further, if $q=m$ \ and \ $\t\ge\max(m,\tt)$, then $v$ has
nonstrict $(\tau,m,\mu,m)$-$N$-property with $\mu=\tau-m$. }
\end{ttt}

It is easy to see, that in formulation of
Theorems~\ref{MSN-D}--\ref{MSN-SL-D} if we take $q=\sigma$, then
$\mu=0$, where $\sigma$ is defined in formulation of the
corresponding Theorems~\ref{MSN}--\ref{MSN-SL}. It means
(see~(\ref{dub3})\,), that Theorems~\ref{MSN-D}--\ref{MSN-SL-D}
include the previous Theorems~\ref{MSN}--\ref{MSN-SL} as
particular case.

\begin{rem}[Some historical remarks]\label{rem-DST}{
It is interesting to note that this Dubovitski\u{\i}
Theorem~A remained almost unnoticed by West mathematicians
for a long time; another proof was given in the recent paper
Bojarski~B. et al.~\cite{BHS}, where they proved also a~version of
this theorem for Holder classes $C^{k,\al+}$ with vanishing
condition~(\ref{uuub}). Further, in ~\cite{HZ} Haj\l{}asz and
Zimmerman replaced the assumption $v\in C^k(\R^n,\R^d)$ of
Theorem~A by the assumption of Sobolev regularity $v\in
W^k_p(\R^n,\R^d)$ with $p>n$ (this is an~analog of DePascale
extension for the~ Morse-Sard, see subsection~\ref{Hist}, cf. with our
assumptions~$kp>n$ or $kp\ge n$ in
theorems~\ref{MSN-S}--\ref{MSN-SL-D}\,).

It is easy to see  that Dubovitski\u{\i}~Theorem~A is a
partial case of Theorem~\ref{MSN-D} of the present paper with
parameters $\tau=n$, $\al=0$ and $q=m+1=d$. Note, that the last
assumption (which also used in~\cite{BHS},~\cite{HZ}\,) simplifies
the proofs very essentially, because automatically one has that
the image $v(E)$ is $\H^q$-$\sigma$-finite. But in general in
theorems~\ref{MSN-D}--\ref{MSN-SL-D} the image $v(E)$ may have
Hausdorff dimension much large than~$q$ for $E\subset \ZZ_{v,m}$ \ with \
$\H^\tau(E)=0$. Nevertheless, the
equality~$\H^\mu(v^{-1}(y)\cap Z_{v,m})=0$ is fulfilled for $q$-almost all~$y\in v(E)$ as required in definition~(\ref{def-dub}\,)

Finally, let us note that the assertions of
Theorems~\ref{MSN-D}--\ref{MSN-SL-D} for the case~$\tau=n$ were
proved in our previous paper~\cite{FKR-MS} and in the papers of
\cite{HKK} by Haj\l{}asz, Korobkov, Kristensen.}
\end{rem}

Without the gradient constraints, the Dubovitski\u{\i}--Fubini analogs of Theorems~\ref{LPT1}--\ref{LPT2} were obtained in our previous paper~\cite{FKR-n}.

\begin{ttt}[\cite{FKR-n}, Sobolev case]\label{LPT1-dub}
{\sl 
Let $\al>0$, \,$1<p<\infty$, \,$\alpha p>n$,  and
$v\in \LK(\R^n,\R^d)$. Suppose that $0<\tau\le n$ and $\tau\ne\tau_*=n-(\alpha-1)p$. Then for every 
$q\in[0,\sigma]$ and for any set $E\subset\R^n$ with $H^\tau(E)=0$
the equality
\begin{equation}\label{dub4-q}
\H^{\mu}(E\cap v^{-1}(y))=0\qquad\mbox{ for \ $\H^q$-a.a.
}y\in\R^d
\end{equation}
holds, where $\mu=\tau\bigl(1-\frac{q}{\sigma}\bigr)$ and the value $\sigma=\sigma(\tau,\al,p)$ is defined in~(\ref{np1}).
}
\end{ttt}

The above Theorem~\ref{LPT1-dub} omits the limiting cases $\alpha p=n$
and $\tau=\tt$. As above, it is  possible to cover these cases as
well using the Lorentz norms.

\begin{ttt}[\cite{FKR-n}, Sobolev--Lorentz case]
\label{LPT2-dub}{\sl  Let  $v:\R^n\to \R^d$ be
a~mapping for which one of the following cases holds:
\begin{itemize}
\item[(i)] \,$v\in W^k_1(\R^n,\R^d)$ for some $k\in\N$, $k\ge n$;

\item[(ii)] \, $v\in \Le^{\al}_{p,1}(\R^n,\R^d)$ \ for  some $\al>0$, \ $p\in (1,\infty)$ with $\al p\ge n$.
\end{itemize}
 Suppose that  $0<\tau\le n$. Then for every 
$q\in[0,\sigma]$ and for any set $E\subset\R^n$ with $H^\tau(E)=0$
the equality~(\ref{dub4-q}) holds with the same~$\mu$ and $\sigma$ defined in~(\ref{np1}) \ (with $\al=k$ and $p=1$ for the case~(i)\,).
}
\end{ttt}

Taking $\tau\ge\tt$, we obtain, in particular, 

\begin{cor}\label{cor-LPT1}
{\sl 
Let $\al>0$, \,$1<p<\infty$, \,$\alpha p>n$,  and
$v\in \LK(\R^n,\R^d)$. Suppose that $0<\tau\le n$ and $\tau>\tau_*=n-(\alpha-1)p$. Then for every 
$q\in[0,\tau]$ and for any set $E\subset\R^n$ with $H^\tau(E)=0$
the equality
\begin{equation}\label{dub4-q1}
\H^{\tau-q}(E\cap v^{-1}(y))=0\qquad\mbox{ for \ $\H^q$-a.a.
}y\in\R^d
\end{equation}
holds. Further, if $v\in\Le^\al_{p,1}(\R^n,\R^d)$ or if $v\in W^k_1(\R^n,\R^d)$, then the same assertion holds under weaker assumptions
$\al p\ge n$ \ (respectively, $k\ge n$\,) and $\tau\ge\tt$. 
}
\end{cor}

\section{Preliminaries}
\label{prel}

\noindent

\noindent By  an {\it $n$--dimensional interval} we mean a closed
cube in $\R^n$ with sides parallel to the coordinate axes. If $Q$
is an $n$--dimensional cubic interval then we write $\ell(Q)$ for
its sidelength.

For a subset $S$ of $\R^n$ we write $\Le^n(S)$ for its outer
Lebesgue measure (sometimes we use the symbol~$\meas S$ for the
same purpose\,). The $m$--dimensional Hausdorff measure is denoted
by $\H^m$ and the $m$--dimensional Hausdorff content by
$\H^{m}_{\infty}$. Recall that for any subset $S$ of $\R^n$ we
have by definition
$$
\H^m (S)=\lim\limits_{t\searrow 0}\H^m_t (S) = \sup_{t >0}
\H^{m}_{t}(S),
$$
where for each $0< t \leq \infty$,
$$
\H^m_t (S)=\inf\left\{ \sum_{i=1}^\infty(\diam S_i)^m\ :\ \diam
S_i\le t,\ \ S \subset\bigcup\limits_{i=1}^\infty S_i \right\}.
$$
It is well known that $\H^n(S)  = \H^n_\infty(S)\sim\Le^n(S)$ for
sets~$S\subset\R^n$.

To simplify the notation, we write $\|f\|_{\LL_p}$ instead of
$\|f\|_{\LL_p(\R^n)}$, etc.

The Sobolev space $\WW^{k}_{p} (\R^n,\R^d)$ is as usual defined as
consisting of those $\R^d$-valued functions $f\in \LL_p(\R^n)$
whose distributional partial derivatives of orders $l\le k$ belong
to $\LL_p(\R^n)$ (for detailed definitions and differentiability
properties of such functions see, e.g., \cite{EG}, \cite{M},
\cite{Ziem}, \cite{Dor}).  We use the norm
$$
\|f\|_{\WW^{k}_{p}}=\|f\|_{\LL_p}+\|\nabla
f\|_{\LL_p}+\dots+\|\nabla^kf\|_{\LL_p},
$$
 and unless otherwise specified all norms on the
spaces $\R^s$ ($s \in \N$) will be the usual euclidean norms.

Working with locally integrable functions, we always assume that
the precise representatives are chosen. If $w\in
L_{1,\loc}(\Omega)$, then the precise representative $w^*$ is
defined for {\em all} $x \in \Omega$ by
\begin{equation*}
\label{lrule}w^*(x)=\left\{\begin{array}{rcl} {\displaystyle
\lim\limits_{r\searrow 0} \dashint_{B(x,r)}{w}(z)\,\dd z}, &
\mbox{ if the limit exists
and is finite,}\\
 0 \qquad\qquad\quad & \; \mbox{ otherwise},
\end{array}\right.
\end{equation*}
where the dashed integral as usual denotes the integral mean,
$$
\dashint_{B(x,r)}{ w}(z) \, \dd
z=\frac{1}{\Le^n(B(x,r))}\int_{B(x,r)}{ w}(z)\,\dd z,
$$
and $B(x,r)=\{y: |y-x|<r\}$ is the open ball of radius $r$
centered at $x$. Henceforth we omit special notation for the
precise representative writing simply $w^* = w$.

If $k<n$, then it is well-known that functions from Sobolev spaces
$\WW^{k}_{p}(\R^n)$ are continuous for $p>\frac{n}k$ and could be
discontinuous for $p\le \p=\frac{n}k$ (see, e.g., \cite{M,Ziem}).
The Sobolev--Lorentz space $\WW^{k}_{\p,1}(\R^n)\subset
\WW^{k}_{\p}(\R^n)$ is a refinement of the corresponding Sobolev
space. Among other things functions that are locally in
$\WW^{k}_{\p,1}$ on $\R^n$ are in particular continuous (see,
e.g., \cite{KK3}\,).

Here we only mentioned the Lorentz space $\LL_{p,1}$, \,$p\ge1$, and in this
case one may rewrite the norm as (see for instance
\cite[Proposition 3.6]{Maly2})
\begin{equation*}\label{lor1}
\|f\|_{L_{p,1}}=
\int\limits_0^{+\infty}\bigl[\Le^n(\{x\in\R^n:|f(x)|>t\})\bigr]^{
\frac1p} \, \dd t.
\end{equation*}
Of course, we have the inequality
\begin{equation}\label{lor1---}
\|f\|_{L_p}\le \|f\|_{L_{p,1}}.\end{equation}

Denote by $\WW^{k}_{p,1}(\R^n)$ the space of all functions $v\in
\WW^k_p(\R^n)$ such that in addition the Lorentz norm~$\|\nabla^k
v\|_{\LL_{p,1}}$ is finite.

By definition put $\|g\|_{L_{p,1}(E)}:=\|1_E\cdot
g\|_{L_{p,1}}$, where $1_E$ is the indicator function of~$E$.  We need  the following analog of the additivity property for the Lorentz norms:
\begin{equation}\label{lor-add}
\sum\limits_i\|f\|^p_{\LL_{p,1}(Q_i)}\le
\|f\|^p_{\LL_{p,1}(\cup_iQ_i)}\qquad \mbox{for any family of disjoint cubes~$Q_i$}
\end{equation}
(see, e.g.,
\cite[Lemma~3.10]{Maly2} or \cite{rom1}\,).

For a function $f\in L_{1,\loc}(\R^n)$ we often use the classical
Hardy--Littlewood maximal function:
$$
\M f(x)=\sup\limits_{r>0}  \dashint_{B(x,r)} \! |f(y)|\,\dd y.
$$

\subsection{On Fubini type theorems  for graphs of continuous functions}

Recall that by usual Fubini theorem, if a set $E\subset\R^2$ has a
zero plane measure, then for $\H^1$-almost all straight lines $L$
parallel to coordinate axes we have $\H^1(L\cap E)=0$. The next
result could be considered as functional Fubini type theorem.

\begin{ttt}[see Theorem~5.3 in~\cite{HKK}]\label{FubN}{\sl
Let $\mu\ge 0$, \,$q>0$, \,and \,$v:\R^n\to\R^d$ \,be a continuous
function. For a set $E\subset\R^n$ define the set function
\begin{equation}\label{dd5}
\Phi(E)=\inf\limits_{E\subset\bigcup_j
D_j}\sum\limits_j\bigl(\diam D_j\bigr)^\mu\bigl[\diam
v(D_j)\bigr]^q,
\end{equation}
where the infimum is taken over all countable families of compact
sets $\{D_j\}_{j\in \N}$ such that $E\subset\bigcup_j D_j$. Then \
$\Phi(\cdot)$ is a countably subadditive and the implication
\begin{equation*}\label{dd6} \Phi(E)=0\ \boldsymbol{\Rightarrow}\
\biggl[\H^\mu\bigl(E\cap v^{-1}(y)\bigr)=0\quad\mbox{for
$\H^q$-almost all }y\in\R^d\biggr]
\end{equation*}holds.}
\end{ttt}

\section{Estimates of the critical values on cubes}\label{appendixII}

In  this section we formulate estimates of the above defined set function $\Phi$ obtained in~\cite[Appendix]{FKR-MS} for subsets of critical set in cubes for different classes of mappings\footnote{The~only technical difference is that in~\cite{FKR-MS} we used the notation $Z'_v=\{x\in\R^n\setminus A_v:\rank\nabla v(x)<m\}$, i.e., there $m-1$ plays the role of the~parameter~$m$ of the~present article.}. 

For all 
the following four subsections  fix $m\in\{0,\dots,n-1\}$ \,and \,$d\ge m$. Take also a positive parameter $q\ge m$  and nonnegative~$\mu\ge0$ required in the definition of the set--function
$\Phi$ in~(\ref{dd5}). 

For a regular (in a sense) 
mapping~$v:\R^n\to\R^d$ denote
$$Z_{v}=\{x\in\R^n\setminus A_v:\rank\nabla v(x)\le m\}.$$ Here $A_v$ means the set of `bad' points, where $v$ is not differentiable or
or which are not Lebesgue points for~$\nabla v$ (of course, $A_v=\emptyset$ if the gradient~$\nabla v$ is a~continuous function\,).  
It is convenient (and sufficient for our purposes) to restrict our attention on the following subset of critical points
\begin{equation}\label{Z'}Z'_v=\{x\in Z_{v}:|\nabla v(x)|\le1\}.
\end{equation}

\subsection{Estimates on cubes for Holder classes of mappings. }\label{HT-es}

Fix \,$k\ge1$, \,$0\le\al\le 1$, \,and
\,$v\in \Cc(\R^n,\R^d)$. By definition of the space
$\Cc,$ there exists a constant~$A\in \R_+$ such that 
\begin{equation}\label{app-hhc1} \mbox{$|\nabla^kv(x)-\nabla^kv(y)|\le  A\cdot |x-y|^\alpha$ \ for all $x,y\in\R^n$.}
\end{equation}

\begin{ttt}[\cite{FKR-MS}]\label{Hold-est-1}{\sl
 Under above assumptions, for any sufficiently small
$n$-dimensional interval $Q\subset\R^n$ the estimate
\begin{equation}
\label{app-6} \Phi(Q\cap Z'_{v})\leq C\,A^{q-m}\,\ell(Q)^{q+\mu+(k+\al-1)(q-m)}
\end{equation}
holds, where the constant $C$ depends on $n,m,k,\alpha,d$ only.}
\end{ttt}

\subsection{Estimates on cubes for Sobolev classes of mappings.  }\label{SC-I}

Fix  \,$k\ge1$, \,$0\le\al< 1$, $1<p<\infty$, \,and
\,$v\in \Lll(\R^n,\R^d)$. In this subsection we consider  the case, when $k+\al>1$ and
\begin{equation}\label{app-sl-1} (k+\al)p> n,
\end{equation}
i.e., when $v$ is  a continuous function (see, e.g., \cite{KK3}\,), but the gradient $\nabla v$ could be discontinuous in general (if $(k+\al-1)p<n$\,).

\begin{ttt}[\cite{FKR-MS}]\label{Sl-est-01}{\sl
 Under above assumptions, there exists a~function $h\in L_{p}(\R^n)$ (depending on~$v$\,) such that the following statements are fulfilled:

 \begin{itemize}
\item[(i)] \,if  $(k+\al-1)p> n$, then the~gradient $\nabla v$ is continuous and uniformly bounded function, and
for any sufficiently small
$n$-dimensional interval $Q\subset\R^n$ the estimate
\begin{equation}
\label{t-sob-2}\Phi(Z'_v\cap Q)\leq C\,\sigma^{q-m}r^{q+\mu+(k+\al-1-\frac{n}p)(q-m)}
\end{equation}
holds, where \begin{equation}
\label{sob-app-s-10} r=\ell(Q),\qquad\sigma=\|h\|_{L_{p}(Q)}.
\end{equation}
and  the constant $C$ depends on $n,m,k,\alpha,d,p$ only.

\item[(ii)] \,if  $(k+\al-1)p< n$, then under additional assumption
 \begin{equation}\label{sl-as}q+\mu>\tt:=n-(k+\al-1)p
\end{equation}for any
$n$-dimensional interval $Q\subset\R^n$ the estimate
\begin{equation}
\label{sob-app-s-9}\Phi(Z'_v\cap Q)\leq C\,\biggl(\,\sigma^q r^{(k+\al-\frac{n}p)q+\mu}+\sigma^{q-m}r^{q+\mu+(k+\al-1-\frac{n}p)(q-m)}\biggr)
\end{equation}
holds with the same $\sigma,r$.
\end{itemize}}
\end{ttt}

\subsection{Estimates on cubes for Sobolev--Lorentz classes of mappings.  }\label{SC-III}

Fix  \,$k\ge1$, \,$0\le\al< 1$, $1<p<\infty$, \,and
\,$v\in \Lll(\R^n,\R^d)$. In this subsection we consider  the case, when $k+\al>1$ and
\begin{equation}\label{app-sl-100} (k+\al)p\ge n,
\end{equation}
i.e., when $v$ is  a continuous function (see, e.g., \cite{KK3}\,), but the gradient $\nabla v$ could be discontinuous in general (if $(k+\al-1)p<n$\,).

\begin{ttt}[\cite{FKR-MS}]\label{Sl-est-3}{\sl
 Under above assumptions, there exists a~function $h\in L_{p,1}(\R^n)$ (depending on~$v$\,) such that the following statements are fulfilled:

 \begin{itemize}
\item[(i)] \,if  $(k+\al-1)p\ge n$, then the~gradient $\nabla v$ is continuous and uniformly bounded function, and
for any sufficiently small
$n$-dimensional interval $Q\subset\R^n$ the estimate~(\ref{t-sob-2}) 
holds with \begin{equation}
\label{sl-app-s-10} r=\ell(Q),\qquad\sigma=\|h\|_{L_{p,1}(Q)}.
\end{equation}

\item[(ii)] \,if  $(k+\al-1)p< n$, then under additional assumption
 \begin{equation}\label{sl-as000}q+\mu\ge\tt:=n-(k+\al-1)p
\end{equation}for any
$n$-dimensional interval $Q\subset\R^n$ the estimate~(\ref{sob-app-s-9}) holds with the same $\sigma,r$ as in~(\ref{sl-app-s-10}).
\end{itemize}}
\end{ttt}

\begin{rem}\label{rem-sl-1}
Formally estimates in Theorem~\ref{Sl-est-3} are the same as in Theorems~\ref{Sl-est-01}, the only difference  is in the definition of~$\sigma$ (using
the~Lorentz norm instead of Lebesgue one\,).
However, Theorem~\ref{Sl-est-3} is `stronger' in  a~sense than the previous Theorems~\ref{Sl-est-01}. Namely, there are
 some important (limiting) cases, which are not covered by Theorem~\ref{Sl-est-01}, but one could still apply the Theorem~\ref{Sl-est-3} for these cases.
It happens for the following values of the parameters:
  \begin{equation}\label{sl-par1}(k+\al)p=n,
\end{equation}
or \begin{equation}\label{sl-par2}(k+\al-1)p=n,
\end{equation}
or
\begin{equation}\label{sl-par3}q+\mu=\tt.
\end{equation}
It means, that the Lorentz norm is a sharper and more accurate tool here than the Lebesgue norm.
\end{rem}

\subsection{Estimates on cubes for Sobolev classes of mappings $W^k_1(\R^n)$, \ $k\ge n$.  }\label{SC-IIII}
In this subsection we consider the limiting case $p=1$ for Sobolev spaces $W^k_1$.
It is well known that functions from the Sobolev space $W^k_1(\R^n,\R^d)$ are continuous if
\begin{equation}
\label{cs-1}k\ge n,
\end{equation}
so we assume this condition below. Fix $k\ge n$  \,and
\,$v\in W^k_1(\R^n,\R^d)$.

\begin{ttt}[\cite{FKR-MS}]\label{Sl-est-4}{\sl
 Under above assumptions, the following statements hold:

 \begin{itemize}
\item[(i)] \,if  $k-1\ge n$, then the gradient $\nabla v$ is continuous and uniformly bounded function, and
for any sufficiently small
$n$-dimensional interval $Q\subset\R^n$ the estimate
\begin{equation}
\label{t-cs-2}\Phi(Z'_v\cap Q)\leq C\,\sigma^{q-m}r^{q+\mu+(k-1-n)(q-m)}
\end{equation}
holds, where again \begin{equation}
\label{cs-app-s-10} r=\ell(Q),\qquad\sigma=\|\nabla^kv\|_{L_{1}(Q)}.
\end{equation}
and  the constant $C$ depends on $n,m,k,d$ only.

\item[(ii)] \,if  $k=n$, then under additional assumption
 \begin{equation}\label{cs-as}q+\mu\ge 1
\end{equation}for any
$n$-dimensional interval $Q\subset\R^n$ the estimate
\begin{equation}
\label{cs-app-s-9}\Phi(Z'_v\cap Q)\leq C\,\biggl(\,\sigma^q r^{\mu}+\sigma^{q-m}r^{\mu+m}\biggr),
\end{equation}
holds with the same $r,\sigma$, and with $C$ depending on $n,m,k,d$ only.
\end{itemize}}
\end{ttt}

\section{Proofs of the main results}\label{proof1}

We have to prove three theorems \ref{MSN-D}--\ref{MSN-SL-D} (because other two theorems \ref{MSN-S}--\ref{MSN-SL} are the partial cases of 
Theorems \ref{MSN-S-D}--\ref{MSN-SL-D} when $q=\sigma$ and $\mu=0$\,). 

For the extremal case $\tau=n$ all these three theorems were proved in~\cite{HKK} and \cite{FKR-MS}, so below we always assume that 
\begin{equation}
\label{pr-f1}0<\t<n.
\end{equation}

Let us first check the assertions about strict $N$-properties. Fixed the corresponding parameters $m\in\{0,1,\dots,n-1\}$, $\mu\ge0$, $q\in(m,\sigma]$, 
 and a~mapping $v:\R^n\to\R^d$ satisfying assumptions of one of the Theorems~\ref{MSN-D}--\ref{MSN-SL-D}. We have to prove that 
\begin{equation}
\label{pr-f2}\H^\mu(E\cap v^{-1}(y))=0\quad\mbox{ for \ $\H^q$-almost all $y\in \R^d$ whenever \ $E\subset\ZZ_{v,m}$ with $\H^\t(E)<\infty$}.
\end{equation}
 
First of all, we will simplify the situation and eliminate some technical difficulties associated with irregular points of mappings from Sobolev classes. Recall, that 
for the Sobolev case the
$m$-critical set is defined as
$$\ZZ_{v,m}=\{x\in\R^n:x\in A_v\mbox{\ \ or \ }x\in\R^n\setminus A_v\  \ \mbox{with}\,\ \rank\nabla v(x)\le m\}.$$
Here $A_v$ means  the set of `bad' points at which either the
function~$v$ is not differentiable or which are not the  Lebesgue
points for~$\nabla v$. Recall that the set $A_v$ is relatively small:
\begin{equation}\label{pr-small1}\H^t(A_v)=0\quad\forall t>\tt:=n-(k+\al-1)p\quad\mbox{if $v\in\Le^{k+\al}_p(\R^n)$ \ (case of Theorem~\ref{MSN-S-D}); }
\end{equation}
\begin{equation}\label{pr-small2}\H^{\tt}(A_v)=0\qquad\mbox{ if \ $v$ is from Theorem~\ref{MSN-SL-D} }.
\end{equation}
In particular, $A_v=\emptyset$ if
$(k+\al-1)p>n$ \ (respectively, if $(k+\al-1)p\ge n$\,).

For the case of Theorem~\ref{MSN-S-D}, take $t\in(\tt,\t)$. Then by Corollary~\ref{cor-LPT1} we have 
\begin{equation}\label{pr-f3}\H^{t-q}(A_v\cap v^{-1}(y))=0\qquad\mbox{ for $\H^q$-almost all $y\in\R^d$}.
\end{equation}
By elementary direct calculation, if $\tt>0$, then 
\begin{equation}
\label{pr-f4}\tt-q<\mu=\t-m-\bigl(k+\al-\frac{n}p+\frac\tau{p}\bigr)(q-m).
\end{equation}
Indeed, by definition of $\tt=n-(k+\al-1)p$, the last inequality is equivalent to
\begin{equation}
\label{pr-f5}(\t-\tt)\bigl(1-\frac{q-m}p\bigr)>0.
\end{equation}
But really by our assumptions 
$$\frac{q-m}p\le\frac{\sigma-m}p=\frac{\t-m}{\t+(k+\al)p-n}<\frac\t\t=1,$$
so (\ref{pr-f4})--(\ref{pr-f5}) is fulfilled. 
From inequality (\ref{pr-f4}) it follows that for $t\in(\tt,\t)$ sufficiently close to $\tt$ we have 
$$t-q<\mu.$$
From this inequality and (\ref{pr-f3}) we obtain
\begin{equation}\label{pr-f6}\H^{\mu}(A_v\cap v^{-1}(y))=0\qquad\mbox{ for $\H^q$-almost all $y\in\R^d$},
\end{equation}
so indeed $A_v$ is negligible in property~(\ref{pr-f2}). 

If $v$ is from Theorem~\ref{MSN-SL-D}, then again Corollary~\ref{cor-LPT1} implies 
\begin{equation}\label{pr-f8}\H^{\tt-q}(A_v\cap v^{-1}(y))=0\qquad\mbox{ for $\H^q$-almost all $y\in\R^d$}.
\end{equation} And by the same calculations we obtain
\begin{equation}
\label{pr-f7}\tt-q\le\mu,
\end{equation}
therefore, the identity~(\ref{pr-f6}) is fulfilled as well and in any case the 'bad' set $A_v$ is negligible in property~(\ref{pr-f2}). 

It means, that in the required property~(\ref{pr-f2}) we could replace the set $\ZZ_{v,m}$ by smaller (regular) set 
$$Z_v=\{x\in\R^n\setminus A_v:\rank\nabla v(x)\le m\}.$$
Moreover, since the countable union of the sets of $\H^\mu$-measure zero has again $\H^\mu$-measure zero, we could replace the set $Z_v$ by the~smaller set
$$Z'_v=\{x\in\R^n\setminus A_v:|\nabla v(x)|\le1\mbox{\ \ and\ \ }\rank\nabla v(x)\le m\},$$
i.e., instead of (\ref{pr-f2}) we need to check only
\begin{equation}
\label{pr-f10}\H^\mu(E\cap v^{-1}(y))=0\quad\mbox{ for \ $\H^q$-almost all $y\in \R^d$ whenever \ $E\subset Z'_v$ with $\H^\t(E)<\infty$}.
\end{equation}
Because of Theorem~\ref{FubN}, for the proof of the last assertion it is sufficient to check, that 
\begin{equation}
\label{pr-f11}\Phi(E)=0\quad\mbox{ whenever \ $E\subset Z'_v$ with $\H^\t(E)<\infty$},
\end{equation}
where the set function $\Phi$ was defined in  Theorem~\ref{FubN}. 

In our previous paper \cite[Appendix]{FKR-MS} we obtained the general estimates 
for~\linebreak $\Phi(Z'_v\cap Q)$, here $Q$ is an~arbitrary $n$-dimensional cube, for all considered cases: Holder, Sobolev (including fractional Sobolev), \,and Sobolev--Lorentz (see their 
formulation in Section~\ref{appendixII} of the present paper).  From these estimates and from the~Holder inequality the required property~(\ref{pr-f11}) follows easily\footnote{Really, the present paper and~\cite{FKR-MS} were written in the same time, so we had in mind 
the~purposes of the present paper when we formulated and proved the~estimates for $\Phi(Z'_v\cap Q)$ in \cite[Appendix]{FKR-MS}. }.

The detailed description of application of these estimates and Holder inequalities was given in \cite{FKR-MS} for the case $\tau=n$. The present case $\tau<n$ is even simpler:
indeed, the most difficult and subtle part in \cite{FKR-MS} was to prove the {\bf strict} $(\tau,q,\mu,m)$-$N$-property for Holder case when $\tau=n$,~---
it requires the application of some generalised Coarea formula, etc. We do not need to touch these difficulties here. The strictness of considered $N$-properties 
for the present case $0<\tau<n$ follows from the following three simple facts:
\begin{equation}\label{pr-uuub}
\mbox{$|\nabla^kv(x)-\nabla^kv(y)|\le\omega(r)\cdot |x-y|^\alpha$
\ \ \quad whenever \ $|x-y|<r$}
\end{equation}
with $\omega(r)\to0$ 
as \,$r\to0$  for $v\in C^{k,\al+}$ or $v\in C^k$ (i.e., $\al=0$\,);
$$\sum\limits_i\|h\|^p_{L_p(Q_i)}\to0\qquad\mbox{ as }\sum\limits_i\ell(Q_i)^\tau\le C, \qquad\sup\limits_{i}\ell(Q_i)\to0,$$ 
for any (fixed) function~$h\in L_p(\R^n)$, where $Q_i$ is a family of nonoverlapping  $n$-dimensional cubes;
$$\sum\limits_i\|h\|^p_{L_{p,1}(Q_i)}\to0\qquad\mbox{ as }\sum\limits_i\ell(Q_i)^\tau\le C, \qquad\sup\limits_{i}\ell(Q_i)\to0,$$ 
for any (fixed) function~$h\in L_{p,1}(\R^n)$, where again $Q_i$ is a family of nonoverlapping  $n$-dimensional cubes (see~(\ref{lor-add})\,).
Since there are no any difficulties  in realisation of these arguments, we omit the details.

The proof of the nonstrict $N$-properties in Theorems~\ref{MSN-D}--\ref{MSN-SL-D} is based on the same estimates with evident 
simplifications in calculations.

\

\noindent {\em Acknowledgment.} M.K. was partially supported by
the Ministry of Education and Science of the Russian Federation
(grant 14.Z50.31.0037) and by the Russian Foundation for Basic Research (Grant 18-01-00649-a).

\noindent Dipartimento di Matematica e Fisica Universit\`a degli
studi della Campania "Luigi Vanvitelli," viale Lincoln 5,
81100, Caserta, Italy\\
e-mail: {\it Adele.FERONE@unicampania.it}

\bigskip

\noindent School of Mathematical Sciences Fudan University,
Shanghai 200433, China, and Voronezh State University,
Universitetskaya pl. 1,
Voronezh, 394018, Russia\\
e-mail: {\it korob@math.nsc.ru}

\bigskip

\noindent Dipartimento di Matematica e Fisica Universit\`a degli
studi della Campania "Luigi Vanvitelli," viale Lincoln 5,
81100, Caserta, Italy\\
e-mail: {\it alba.roviello@unicampania.it}

\end{document}